\documentclass[12pt,reqno,a4wide]{amsart}

\usepackage{tikz} 
\usepackage{calrsfs}
\usepackage{mathrsfs}

\usepackage{array}
\usepackage{hyperref}

\usetikzlibrary{decorations.pathreplacing}
\usetikzlibrary{fit}

\usepackage{pgfplots}

\allowdisplaybreaks

\begin{document}
\bibliographystyle{plain}

%
%

	\title 
	{Dispersed Dyck paths  revisited}

	\author[H. Prodinger ]{Helmut Prodinger }
	\address{Department of Mathematics, University of Stellenbosch 7602, Stellenbosch, South Africa
	and
NITheCS (National Institute for
Theoretical and Computational Sciences), South Africa.}
	\email{hproding@sun.ac.za}

	\keywords {Dispersed Dyck paths,  prefix, ascent, descent, valley kernel method}
	
	\begin{abstract}
Dispersed Dyck paths are Dyck paths, with possible flat steps on level 0. We revisit and augment questions about them
from the Encyclopedia of Integer Sequences, in a systematic way that uses generating functions and the kernel method.
	\end{abstract}
	
	\subjclass[2010]{05A15}

\maketitle

\section{Introduction}


Dispersed Dyck paths consist of up-steps $U=(1,1)$ and down-steps $D=(1,-1)$, never go below the $x$-axis, and can have horizontal steps on the $x$-axis,
$(n,0)\to(n+1,0)$. According to the OEIS \cite{OEIS}, it seems that they were invented/suggested by Emeric Deutsch. 
All sequences that are explicitly cited are from \cite{OEIS}, using the local identifiers $A\#\#\#\#\#\#$.
 
 As in various examples discussed in the past \cite{Prodinger-kernel, Prodinger-prefix, garden} we will use the kernel method to set up appropriate
 generating functions. There are three (or four) bivariate (or trivariate) generating functions $F(u), G(u), H(u)$ (also depending on $z$) related to the nature of the last step of the 
 prefix of a dispersed Dyck path. To solve the system, a `bad' factor $u-r_2$ must be divided out from numerator and denominator, after which one can plug in $u=0$
 and identify the unknown quantities $f_1=F(0)$, etc. 
 
 Various questions (from the encyclopedia of integer sequences  \cite{OEIS}) about such paths will be revisited; our method is fairly automatic and purely refers to generating functions. To be more specific, these include 1-ascents, 1-descents, valleys on level 0, occurrences of $UUDD$.
 Also, we obtain expression for paths that do not have to go back to the $x$-axis, rather finish at a level $j$, or, more generally, on \emph{any} level. This is either achieved by looking at the coefficient of $u^j$, say, or, setting $u=1$. The generating functions of interest usually have 3 variables: $z$ for the length (number of steps), $u$ for the final height, and $t$ for an additional parameter of interest.

\section{Counting 1-ascents} 

A 1-ascent is an ascent consisting of exactly 1 up step. The paper \cite{theis} contains some analysis, but without generating functions. The paper \cite{hackl.AofA} discusses
$d$-ascents for $d\ge2$ as well, but the quadratic equations that are so common in the context of Dyck paths, are then of higher order, and the results are consequently of an
asymptotic nature.
 
 We distinguish 3 states, together with the current level. Down-steps are unproblematic, but when after them an up-step arrives, it might be the only one or further up-steps follow.
 A graph describes all the possible scenarios, as in \cite{Prodinger-prefix}.  It has 3 layers of states, and
  sequences $f_i$, $i\ge1$,  $g_i$, $h_i$, $i\ge0$, in that order. These quantities all depend on the variable $z$ and
 describe generating functions of paths leading to a particular state.  Since on level 0, flat (horizontal) steps are also allowed, the quantities $f_0$, $g_0$, and $h_0$ are somewhat special and will be treated as parameters. Compare \cite{Prodinger-JIS} for this technique. We will deal with trivariate generating functions
 \begin{equation*}
F(u)=\sum_{i\ge1}f_iu^{i-1},\ G(u)=\sum_{i\ge1}g_iu^{i-1},\ H(u)=\sum_{i\ge1}h_iu^{i-1},
 \end{equation*}
 the other variables $z$ (counting the length) and $t$ (counting the 1-ascents) are not explicitly mentioned.
  \begin{figure}[h]

 	\begin{center}
 		\begin{tikzpicture}[scale=1.8,main node/.style={circle,draw,font=\Large\bfseries}]

 			\foreach \x in {0,1,2,3,4,5,6,7,8}
 			{
 				\draw (\x,0) circle (0.05cm);
 				\fill (\x,0) circle (0.05cm);
 				\draw (\x,-1) circle (0.05cm);
 				\fill (\x,-1) circle (0.05cm);
 			}

 			\fill (0,0) circle (0.08cm);

 			\foreach \x in {0,1,2,3,4,5,6,7}
 			{
 				\draw[ thick, -latex] (\x+1,0) to  (\x,0);	
\draw[ thick,  -latex] (\x,0) to  (\x+1,-1);	
\draw[ thick,  -latex] (\x,-1)   to  (\x+1,-2);
\draw[ thick,  -latex] (\x,-2) to  (\x+1,-2);
 			}

 		\foreach \x in {1,2,3,4,5,6,7}
 	{
 		\draw[ thick, -latex] (\x+1,-2) to  (\x,0);	
 	}

 			\draw [ thick,cyan, cyan, -latex](0,-0) .. controls (-0.5,-0.5)  and (-0.5,0.5) .. (0,0);


 			\foreach \x in {0,1,2,3,4,5,6,7,8}
{
	\draw (\x,0) circle (0.05cm);
	\fill (\x,0) circle (0.05cm);
	\draw (\x,-1) circle (0.05cm);
	\fill (\x,-1) circle (0.05cm);
}

\foreach \x in {1,2,3,4,5,6,7,8}
{
	\draw (\x,0) circle (0.05cm);
	\fill (\x,0) circle (0.05cm);
}

	\foreach \x in {0,1,2,3,4,5,6,7,8}
{
	\draw (\x,0) circle (0.05cm);
	\fill (\x,0) circle (0.05cm);
	\draw (\x,-2) circle (0.05cm);
	\fill (\x,-2) circle (0.05cm);
}

\fill (0,0) circle (0.08cm);

\foreach \x in {0,1,2,3,4,5,6,7}
{ 				
	\draw  [thick,red, -latex] (\x+1,-1)[out=110, in=-30] to  (\x,0);	
}

 	\draw [ thick,cyan, cyan, -latex](0,-0) .. controls (-0.5,-0.5)  and (-0.5,0.5) .. (0,0); 		
 			
 		\end{tikzpicture}
 	\caption{Three layers of states, labelled $f,g,h$, in that order. The state $f_0$ is responsible for dispersed Dyck paths, and all others to prefixes of them.}
 	\end{center}
 \end{figure}
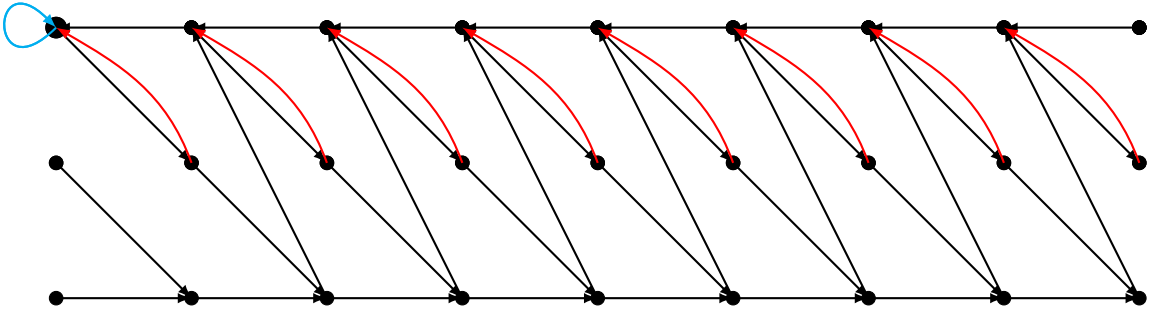
The recursions can be read off, by considering the last step made,
\begin{align*}
f_i&=zf_{i+1}+ztg_{i+1}+zh_{i+1},\ i\ge 1,\ f_0=1+zf_0+zf_1+ztg_1, \\
g_{i+1}&=zf_i,\ i\ge1,\ g_{1}=zf_0,\\
h_{i+1}&=zg_i+zh_i,\ i\ge0.
\end{align*}
Translating these into the trivariate generating functions, we get
\begin{align*}
\sum_{i\ge1}u^{i}f_i&=\sum_{i\ge1}u^{i}zf_{i+1}+\sum_{i\ge1}u^{i}ztg_{i+1}+\sum_{i\ge1}u^{i}zh_{i+1},\\
uF(u)&=zF(u)-zf_1+ztG(u)-ztg_1+zH(u),
\end{align*}
\begin{equation*}
G(u)-g_1=zuF(u),\quad H(u)=zuG(u)+zuH(u).
\end{equation*}
The system can be solved, but still depends on the initial values $f_1, g_1, h_1$:
\begin{align*}
F(u)&={\frac {z  ( -f_1+zuf_1+zug_1  ) }{-{z}^{3}{u}
		^{2}+u-z{u}^{2}-z+{z}^{2}u-{z}^{2}tu+{z}^{3}t{u}^{2}}},\\
G(u)&={\frac {-{z}^{2}uf_1-{z}^{2}tug_1-zg_1+ug_1+{z}^{3
		}{u}^{2}f_1+{z}^{3}t{u}^{2}g_1+{z}^{2}ug_1-z{u}^{2}g_1}{-{z}^{3}{u}^{2}+u-z{u}^{2}-z+{z}^{2}u-{z}^{2}tu+{z}^{3}t{u}^{			2}}},\\
H(u)&=-{\frac {zu  ( {z}^{2}uf_1+{z}^{2}tug_1+zg_1-ug_1  ) }{-{z}^{3}{u}^{2}+u-z{u}^{2}-z+{z}^{2}u-{z}^{2}tu+{z}^{3}t
		{u}^{2}}}.
\end{align*}
The naive approach would be to plug in $u=0$ and identity the initial values, but this doesn't work, and requires some preparations.
The denominator will be factored,
\begin{equation*}
-{z}^{3}{u}
^{2}+u-z{u}^{2}-z+{z}^{2}u-{z}^{2}tu+{z}^{3}t{u}^{2}=z(-z^2-1+z^2t)(u-r_1)(u-r_2)
\end{equation*}
with
\begin{equation*}
r_1=\frac{1+z^2(1-t)+W}{2z(1+z^2(1-t))},\quad
r_2=\frac{1+z^2(1-t)-W}{2z(1+z^2(1-t))}
\end{equation*}
and the abbreviation
\begin{equation*}
	W:=\sqrt{1-2(t+1)z^2-(t+3)(1-t)z^4}.
\end{equation*}
Dividing out the factor $u-r_2$  from numerator and denominator (`kernel method'), we find
\begin{align*}
F(u)&={\frac {{z}^{2}  ( f_1+g_1  ) }{  ( -{z}^{2}-1+
		{z}^{2}t  )   ( r_2 z+zu-1  ) }},\\*
G(u)&={\frac {r_2 {z}^{3}f_1+r_2 {z}^{3}tg_1-r_2
		 zg_1-{z}^{2}f_1-{z}^{2}tg_1+g_1+{z}^{2}g_1+
		u{z}^{3}f_1+u{z}^{3}tg_1-zug_1}{  ( -{z}^{2}-1+{z}^
		{2}t  )   ( r_2 z+zu-1  ) }},\\
H(u)&=-{\frac {z  ( r_2 {z}^{2}f_1+r_2 {z}^{2}tg_1
		-r_2 g_1+{z}^{2}uf_1+{z}^{2}tug_1+zg_1-ug_1  ) }{  ( -{z}^{2}-1+{z}^{2}t  )   ( r_2		 z+zu-1  ) }},
\end{align*}
and now $u=0$ is possible, with
\begin{align*}
f_1&={\frac {{z}^{2}  ( f_1+g_1  ) }{  ( -{z}^{2}-1+
		{z}^{2}t  )   ( r_2 z-1  ) }},\quad
g_1=zf_0,\quad h_1=0.
\end{align*}
But we have
\begin{equation*}
f_0=1+zf_0+zf_1+ztg_1,
\end{equation*}
and both, $f_1$ and $g_1$ can be expressed in terms of $f_0$:
\begin{equation*}
f_1={\frac {{z}^{2}g_1}{1-{z}^{2}t-r_2 z-r_2 {z}^{3}+r_2 {z}^{3}t}}.
\end{equation*}
The ultimate solution is now
\begin{equation*}
f_0= {\frac {-1+2 z-{z}^{2}+{z}^{2}t+W}{2z  ( {z}^{2}+1
		-2 z-{z}^{3}+{z}^{3}t-{z}^{2}t  ) }},
\end{equation*}
with the series expansion
\begin{equation*}
f_0=1+z+ ( t+1 ) {z}^{2}+ ( 2 t+1 ) {z}^{3}+
( {t}^{2}+3 t+2 ) {z}^{4}+ ( 3 {t}^{2}+4 t+3) {z}^{5}+ ( {t}^{3}+6 {t}^{2}+8 t+5 ) {z}^{6}+\cdots.
\end{equation*}
For $t=0$ we find the enumeration of dispersed Dyck paths without 1-ascents (A191385)
\begin{align*}
f_0\big|_{t=0}&=\frac{-1+2z-z^2+\sqrt{1-2z^2-3z^4}}{2z(1-2z+z^2-z^3)}\\
&=1+z+z^2+z^3+2z^4+3z^5+5z^6+7z^7+12z^8+18z^9+31z^{10}+47z^{11}+81z^{12} +\cdots.
\end{align*}
For $t=1$, we find
\begin{equation*}
f_0\big|_{t=1}=-\frac1{2z}+\frac{\sqrt{1-4z^2}}{2z(1-2z)}=\sum_{n\ge0}\binom{2n}{n}z^{2n}+\sum_{n\ge0}\binom{2n+1}{n}z^{2n+1},
\end{equation*}
which is the enumeration of dispersed Dyck paths of length $n$: $\binom{n}{\lfloor n/2\rfloor}$.
The number of 1-ascents in paths of length 5 can be deduced from $( 3 {t}^{2}+4 t+3)$: it is $3\cdot2+4$, and in general, we must differentate 
$f_0$ w.r.t. $t$, followed by $t=1$.
\begin{align*}
	\frac{\partial f_0}{\partial z}\bigg|_{t=1}&=\frac{z^2(1-4z^2+\sqrt{1-4z^2})}{2(1-2z)(1-4z^2)}\\
	&=\frac{z^2}{2(1-2z)}+\frac{z^2}{2(1-4z^2)^{3/2}}+\frac{z^3}{(1-4z^2)^{3/2}}\\
	&=\sum_{n\ge2}2^{n-3}z^n+\sum_{n\ge1} \frac{(2n-1)!}{2(n-1)!(n-1)!}z^{2n}+\sum_{n\ge1} \frac{(2n-1)!}{(n-1)!(n-1)!}z^{2n+1}.
\end{align*}
The coefficients of $z^n$ form the sequene A045621.
Furthermore, we get
\begin{align*}
F(u)&=\frac{z^2r_2}{1-zr_2-zu}f_0,\\*
G(u)&=\frac{z(1-zu)(1-zr_2)}{1-zr_2-zu}f_0,\\*
H(u)&=\frac{uz^2(1-zr_2)}{1-zr_2-zu}f_0.
\end{align*}
From this, it is easy to find $f_j, g_j, h_j$ since we only have to expand $1/(1-zr_2-zu)$ in powers of $u$, which is basically a geometric series. We will not write this out.
However, we will sum \emph{all} these quantities, as it describes all dispersed Dyck paths with open end (partial paths, sometimes called meander).
We have to be careful and, when a path ends in a state from the second layer (`g'), the last up step is a 1-ascent, and we need to attach an extra factor $t$.
The resulting formula is surprisingly simple:
\begin{equation*}
\frac{\partial}{\partial t}\Big(f_0+F(1)+tG(1)+H(1)\Big) \bigg|_{t=1}=\frac{z(1-z)^2}{(1-2z)^2}=z+\sum_{n\ge2}(n+2)2^{n-3}z^n.
\end{equation*}

\section{Counting 1-descents}

In this section, the number of 1-descents will be counted. The definition is similar; a down-step rendered by up-steps or standing at the end of the  dispersed Dyck path.
The same letters as in the previous section will be used, but now with a different meaning. The figure is again self-explanatory:
\begin{figure}[h]

	\begin{center}
		\begin{tikzpicture}[scale=1.8,main node/.style={circle,draw,font=\Large\bfseries}]

			\foreach \x in {0,1,2,3,4,5,6,7,8}
			{
				\draw (\x,0) circle (0.05cm);
				\fill (\x,0) circle (0.05cm);
				\draw (\x,-1) circle (0.05cm);
				\fill (\x,-1) circle (0.05cm);
			}

			\fill (0,0) circle (0.08cm);

			\foreach \x in {0,1,2,3,4,5,6,7}
			{
				\draw[ thick, latex-] (\x+1,0) to  (\x,0);	
				\draw[ thick,  -latex] (\x+1,0) to  (\x,-1);	
				\draw[ thick,  -latex] (\x+1,-1)   to  (\x,-2);
				\draw[ thick,  -latex] (\x+1,-2) to  (\x,-2);
			}

			\foreach \x in {0,1,2,3,4,5,6,7}
			{
				\draw[ thick, -latex] (\x,-2) to  (\x+1,0);	
			}

			\draw [ thick,cyan, cyan, -latex](0,-0) .. controls (-0.5,-0.5)  and (-0.5,0.5) .. (0,0);

			
			\foreach \x in {0,1,2,3,4,5,6,7,8}
			{
				\draw (\x,0) circle (0.05cm);
				\fill (\x,0) circle (0.05cm);
				\draw (\x,-1) circle (0.05cm);
				\fill (\x,-1) circle (0.05cm);
			}
			
			\foreach \x in {1,2,3,4,5,6,7,8}
			{
				\draw (\x,0) circle (0.05cm);
				\fill (\x,0) circle (0.05cm);
			}
			
			\foreach \x in {0,1,2,3,4,5,6,7,8}
			{
				\draw (\x,0) circle (0.05cm);
				\fill (\x,0) circle (0.05cm);
				\draw (\x,-2) circle (0.05cm);
				\fill (\x,-2) circle (0.05cm);
			}
			
			\fill (0,0) circle (0.08cm);
			
			\foreach \x in {0,1,2,3,4,5,6,7}
			{ 				
				\draw  [thick,red, -latex] (\x,-1)[out=80, in=-160] to  (\x+1,0);	
			}
			\draw [ thick,cyan, cyan, -latex](0,-0) .. controls (-0.5,-0.5)  and (-0.5,0.5) .. (0,0); 		
			\draw [ thick,cyan, cyan, -latex](0,-1) .. controls (-0.5,-0.5-1)  and (-0.5,0.5-1) .. (0,-1); 	
			\draw [ thick,cyan, cyan, -latex](0,-2) .. controls (-0.5,-0.5-2)  and (-0.5,0.5-2) .. (0,-2); 	
		\end{tikzpicture}
		\caption{Three layers of states, labelled $f,g,h$, in that order.}
	\end{center}
\end{figure}
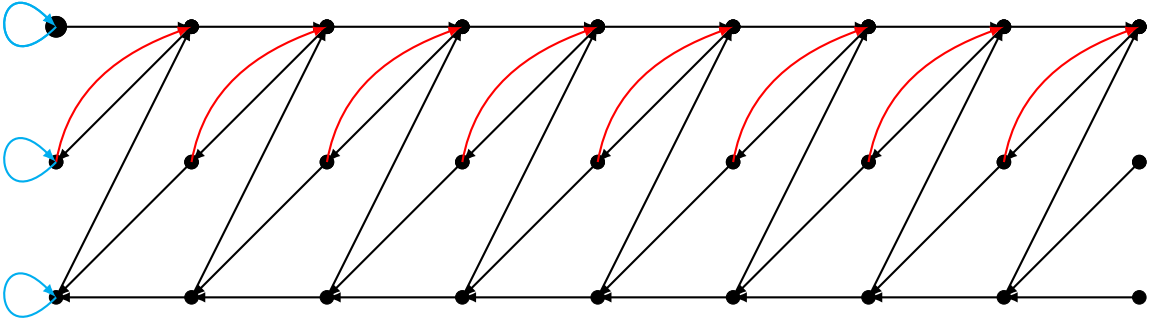
Here are the relevant recursions:
\begin{align*}
f_{i+1}&=zf_i+ztg_i+zh_i,\ i\ge1,\ f_0=\frac1{1-z},\\*
g_i&=zf_{i+1},\ i\ge1,\ g_0=zg_0+zf_1,\ g_0=\frac{z}{1-z}f_1\\*
h_i&=zg_{i+1}+zh_{i+1},\ i\ge1,\ h_0=zh_0+zg_1+zh_1,\ h_0=\frac{z}{1-z}(g_1+h_1),
\end{align*}
and
\begin{align*}
f_{1}&=zf_0+ztg_0+zh_0=\frac z{1-z}+\frac{z^2t}{1-z}f_1+\frac{z^2}{1-z}(g_1+h_1)\\
&=\frac{z(1+zg_1+zh_1)}{1-z-z^2t}.
\end{align*}
So the kernel method will give us $g_1$ and $h_1$, and $f_1$ follows. The recursions will be 
translated into trivariate generating functions, as before:
\begin{align*}
	\sum_{i\ge1}u^if_{i+1}&=\sum_{i\ge1}u^izf_i+\sum_{i\ge1}u^iztg_i+\sum_{i\ge1}u^izh_i,\\
	F(u)-f_1&=uzF(u)+uztG(u)+uzH(u);
\end{align*}
\begin{align*}
	\sum_{i\ge1}u^ig_i&=\sum_{i\ge1}u^izf_{i+1},\\
	uG(u)&=zF(u)-zf_1;
\end{align*}
\begin{align*}
\sum_{i\ge1}u^ih_i&=\sum_{i\ge1}u^izg_{i+1}+\sum_{i\ge1}u^izh_{i+1},\\
uH(u)&=zG(u)-zg_1+zH(u)-zh_1.
\end{align*}
We do not write the solutions for $F(u),G(u),H(u)$, only the relevant denominator
\begin{equation*}
z^3-u+u^2z+uz^2t+z-uz^2-z^3t=z(u-r_1)(u-r_2)
\end{equation*}
with
\begin{equation*}
r_1=\frac{1+z^2-z^2t+W}{2z},\quad r_2=\frac{1+z^2-z^2t-W}{2z},
\end{equation*}
and
\begin{equation*}
W=\sqrt{1-2z^2t-2z^2+z^4t^2+2z^4t-3z^4}.
\end{equation*}
Then we divide out the factor $u-r_2$ in the usual way, and get, after setting $u=0$,
\begin{equation*}
g_1=\frac{zr_2}{1+z^2-z^2t}f_1 ,\quad g_1+h_1=\frac{r_2-z}{z}f_1,
\end{equation*}
\begin{equation*}
f_1=\frac{z}{1-z+z^2-z^2t-zr_2}.
\end{equation*}
This provides also the values $f_0,g_0,h_0$, and
\begin{equation*}
f_0+tg_0+h_0=	\frac{1-r_2}{1-2z+z^2-z^3-z^2t+z^3t},
\end{equation*}
which is the generating function of dispersed Dyck paths returning to the 0-level. Again, to count the contributions of the 1-descents, we compute
\begin{align*}
\frac{\partial}{\partial z}(f_0+tg_0+h_0)\Big|_{t=1}&=\frac{z^2(1-4z^2+\sqrt{1-4z^2})}{2(1-2z)(1-4z^2)}\\
&=\frac{z^2}{2(1-2z)}+\frac{z^2}{2(1-4z^2)^{3/2}}+\frac{z^3}{(1-4z^2)^{3/2}}
\end{align*}
as before. This is natural, as reading from right to left turns 1-upsteps to 1-downsteps, and vice versa. Considering \emph{all} paths, regardless where
they end, 
\begin{equation*}
\frac{\partial}{\partial z}(f_0+F(1)+tg_0+tG(1)+h_0+H(1))\Big|_{t=1}=\frac{z^2}{2(1-2z)^2}+\frac{z^2\sqrt{1-4z^2}}{2(1-2z)^2},
\end{equation*}
which is surprisingly simple, with a simple series expansion,
\begin{equation*}
\sum_{n\ge2}(n-1)2^{n-3}z^n+\sum_{n\ge1}\frac{2(2n-2)!}{(n-1)!(n-1)!}z^{2n}+\sum_{n\ge1}\frac{(2n-2)!(4n-3)}{2(n-1)!(n-1)!}z^{2n+1}.
\end{equation*}

\section{Counting valleys on level 0}

We do now exactly what the title of the section says;
\begin{figure}[h]

	\begin{center}
		\begin{tikzpicture}[scale=1.8,main node/.style={circle,draw,font=\Large\bfseries}]

			\fill (0,0) circle (0.08cm);

			\foreach \x in {0,1,2,3,4,5,6,7}
			{
				\draw[ thick, latex-] (\x+1,0) to  (\x,0);	
								
			}			
			\foreach \x in {1,2,3,4,5,6,7}
{
	\draw[ thick, latex-] (\x,0) [ in =130, out=40]to  (\x+1,0);	
}

							\draw[ thick, -latex] (1,0) to  (0,-1);

			\draw [ thick,cyan, -latex](0,-0) .. controls (-0.5,-0.5)  and (-0.5,0.5) .. (0,0); 
			\draw [ thick,cyan,  latex-](0,-0) to (0,-1);


			\foreach \x in {1,2,3,4,5,6,7,8}
			{
				\draw (\x,0) circle (0.05cm);
				\fill (\x,0) circle (0.05cm);
			}
			
			\foreach \x in {0,1,2,3,4,5,6,7,8}
			{
				\draw (\x,0) circle (0.05cm);
				\fill (\x,0) circle (0.05cm);
			}
			\foreach \x in {0}
{
	\draw (\x,-1) circle (0.05cm);
	\fill (\x,-1) circle (0.05cm);
}
			
			\fill (0,0) circle (0.08cm);
			
			\foreach \x in {0}
			{ 				
				\draw  [thick,red, -latex] (\x,-1)[out=80, in=-160] to  (\x+1,0);	
			}
		 	\draw [ thick,cyan, -latex] (0,-0) .. controls (-0.5,-0.5)  and (-0.5,0.5) .. (0,0); 		
		\end{tikzpicture}
		\caption{Two layers of states, labelled $f,g$, in that order.}
	\end{center}
\end{figure}
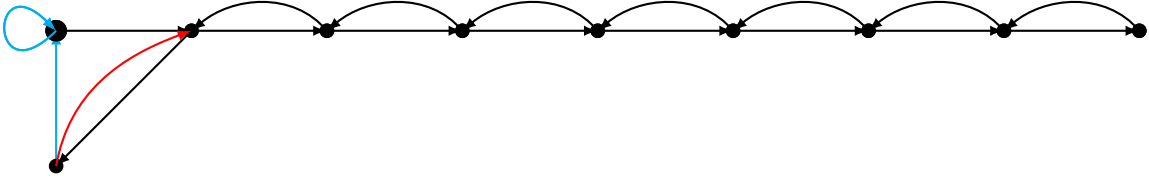

Here are the usual recursions,
\begin{equation*}
f_0=1+zf_0+zg_0,\ g_0=zf_1,\ f_1=ztg_0+zf_0+zf_2,\ f_i=zf_{i-1}+zf_{i+1},\ i\ge2.
\end{equation*}
\begin{equation*}
f_0=\frac{1}{1-z}+\frac{z}{1-z}g_0=\frac{1}{1-z}+\frac{z^2}{1-z}f_1,\quad tf_0-t-ztf_0=ztg_0
\end{equation*}
\begin{align*}
		\sum_{i\ge2}u^{i-1}f_i&=\sum_{i\ge2}u^{i-1}zf_{i-1}+\sum_{i\ge2}u^{i-1}zf_{i+1}\\
		F(u)-f_1&=uzF(u)+\frac zu(F(u)-f_1) -zf_2,\\
		F(u)&=uzF(u)+\frac zu(F(u)-f_1) +ztg_0+zf_0,\\
				F(u)&=uzF(u)+\frac zu(F(u)-f_1) +z^2tf_1+zf_0.
							 \end{align*}
						 We have $f_0=\frac{1+z^2f_1}{1-z}$, and from the kernel method, after dividing out
						 \begin{equation*}
u-r_2,\quad \text{with}\quad r_2:=\frac{1-\sqrt{1-4z^2}}{2z},
						 \end{equation*}
\begin{equation*}
f_1=\frac {z}{1-z-{z}^{2}t+{z}^{3}(t-1)+z(z-1)r_2}.
\end{equation*}
Then,
\begin{align*}
f_1\big|_{t=0}&=\frac{2}{2-3z+z\sqrt{1-z^2}}\\
&=1+z+z^2+2z^3+3z^4+5z^5+8z^6+14z^7+23z^8+41z^9+69z^{10}+125z^{11}+\cdots
\end{align*}
which is sequence A191388 (no valleys on level 0). Furthermore
\begin{equation*}
F(u)=\frac{f_1}{1-ur_2} \Longrightarrow [u^{j-1}]F(u)=f_j=f_1r_2^{j-1}.
\end{equation*}
Finally,
\begin{align*}
\frac{\partial f_0}{\partial t}\bigg|_{t=1}&=\frac{1-3z^2+(z^2-1)\sqrt{1-4z^2}}{2z(1-2z)}\\
&=z^5+2z^6+7z^7+14z^8+37z^9+74z^{10}+176z^{11}+352z^{12}+\cdots,
\end{align*}
which is sequence A191389 (number of  valleys on level 0).  By setting $u=1$ in $g_0+F(u)$, we can address the number of such
partial dispersed Dyck paths according to no valleys resp.\ number of valleys. Since these formul\ae\ are easy to obtain and not
too attractive, we do not display them here.

\section{Counting occurrences of $UUDD$}

\begin{figure}[h]

	\begin{center}
		\begin{tikzpicture}[scale=1.8,main node/.style={circle,draw,font=\Large\bfseries}]

			\fill (0,0) circle (0.08cm);

\foreach \x in {0,1,2,3,4,5,6,7}
			{
				\draw[ thick, -latex] (\x+1,0) to  (\x,0);	
				
			}			

\foreach \x in {0,1,2,3,4,5,6,7}
{
	\draw[ thick, latex-] (\x+1,-1) to  (\x,0);	
	
}			
	\foreach \x in {1,2,3,4,5,6,7}
	{
		\draw[ thick, latex-] (\x+1,-2) to  (\x,-1);	
		
	}			

\foreach \x in {2,3,4,5,6,7}
{
	\draw[ thick, latex-] (\x+1,-2) to  (\x,-2);	
	
}	

\foreach \x in {1,2,3,4,5,6,7}
{
	\draw[ thick, -latex] (\x,-3) to  (\x+1,-1);	
	
}

\foreach \x in {0,1,2,3,4,5,6,7}
{
	\draw[ thick, -latex] (\x+1,-1) [out=100, in=-20]to  (\x,0);	
	
}
	\foreach \x in {1,2,3,4,5,6,7}
{
	\draw[ thick, -latex] (\x+1,-2) to  (\x,-3);	
	
}			\fill (0,0) circle (0.08cm);
			
	\foreach \x in {0,1,2,3,4,5,6,7}
{
	\draw[ thick,red, -latex] (\x+1,-3) to  (\x,0);	
	
}			\fill (0,0) circle (0.08cm);
			\draw [ thick,cyan, -latex] (0,-0) .. controls (-0.5,-0.5)  and (-0.5,0.5) .. (0,0); 
			
			\foreach \x in {1,2,3,4,5,6,7,8}
			{
				\draw (\x,0) circle (0.05cm);
				\fill (\x,0) circle (0.05cm);
			}
		\foreach \x in {1,2,3,4,5,6,7,8}
		{
			\draw (\x,-1) circle (0.05cm);
			\fill (\x,-1) circle (0.05cm);
		}
		\foreach \x in {2,3,4,5,6,7,8}
{
	\draw (\x,-2) circle (0.05cm);
	\fill (\x,-2) circle (0.05cm);
}
		\foreach \x in {1,2,3,4,5,6,7,8}
{
	\draw (\x,-3) circle (0.05cm);
	\fill (\x,-3) circle (0.05cm);
}

		\end{tikzpicture}
		\caption{Four layers of states, labelled $f,g,h,k$, in that order.}
	\end{center}
\end{figure}
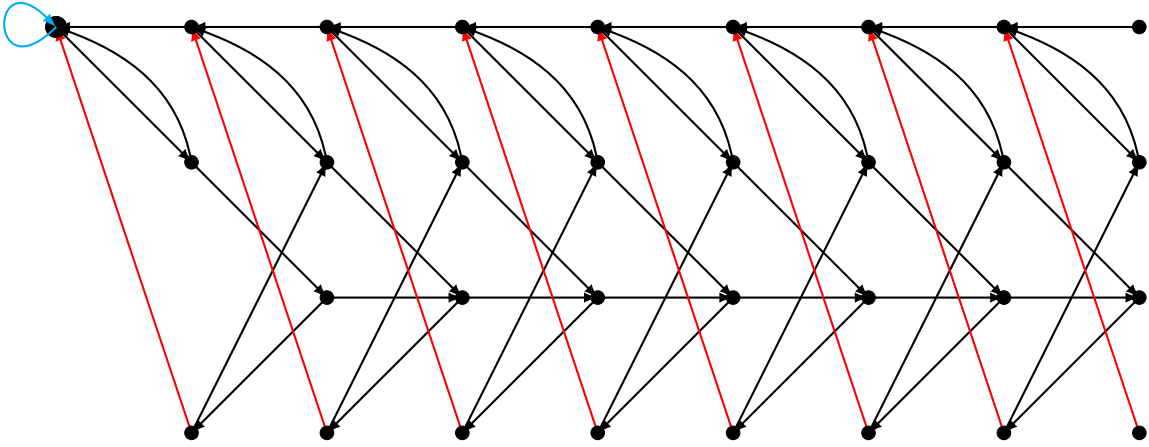
The recursions are
\begin{align*}
f_i&=zf_{i+1}+zg_{i+1}+ztk_{i+1},\ i\ge1,\ f_0=1+zf_0+zf_{1}+zg_{1}+ztk_{1},\\
g_{i+1}&=zf_i+zk_i,\ i\ge1,\ g_1=zf_0\\
h_{i+1}&=zg_i+zh_i,\ i\ge2,\ h_2=zg_1,\\
k_{i}&=zh_{i+1},\ i\ge1,\ k_1=zh_2=z^2g_1.
\end{align*}
The following generating functions will be used:
\begin{align*}
	F(u)=\sum_{i\ge1}f_iu^{i-1},\ G(u)=\sum_{i\ge1}g_iu^{i-1},\ H(u)=\sum_{i\ge2}h_iu^{i-2},\ K(u)=\sum_{i\ge1}k_iu^{i-1}.
\end{align*}
Summing and simplifying,
\begin{align*}
	\sum_{i\ge1}u^if_i&=\sum_{i\ge1}u^izf_{i+1}+\sum_{i\ge1}u^izg_{i+1}+\sum_{i\ge1}u^iztk_{i+1},\\
	uF(u)&=zF(u)-zf_1+zG(u)-zg_1+ztK(u)-ztk_1,\\
	\sum_{i\ge1}u^ig_{i+1}&=\sum_{i\ge1}u^izf_i+\sum_{i\ge1}u^izk_i,\\
	G(u)-g_1&=zuF(u)+zuK(u),\\
	\sum_{i\ge2}u^{i-1}h_{i+1}&=\sum_{i\ge2}u^{i-1}zg_i+\sum_{i\ge2}u^{i-1}zh_i,\\
	H(u)-h_2&=zG(u)-zg_1+zuH(u), \ H(u)=\frac{z}{1-zu}G(u),\\
	\sum_{i\ge1}u^{i-1}k_{i}&=\sum_{i\ge1}u^{i-1}zh_{i+1},\\
K(u)&=zH(u).
\end{align*}
It is beneficial to reduce the system to just one equation, for $F(u)$. Note also that
\begin{equation*}
f_1= -\frac{1+zf_0+z^2f_0+z^4tf_0-f_0}{z},\ g_1=zf_0,\ h_2=z^2f_0,\  k_1=z^3f_0.
\end{equation*}
Then
\begin{equation*}
F(u)=-{\frac {-{z}^{3}u-{z}^{4}f_0u-zu-{z}^{2}f_0u+zuf_0
		+1+zf_0+{z}^{2}f_0+{z}^{4}tf_0-f_0}{-{z}^{4}u+z+{z		}^{4}tu+z{u}^{2}-u}}.
\end{equation*}
Dividing out, as part of the kernel method,
\begin{equation*}
u-r_2,\quad\text{with}\quad r_2= {\frac {1+{z}^{4}-{z}^{4}t-\sqrt {{z}^{8}-2 {z}^{8}t+2 {z}^{4			}+{z}^{8}{t}^{2}-2 {z}^{4}t+1-4 {z}^{2}}}{2z}}
\end{equation*}
and setting $u=0$ leads to
\begin{equation*}
-\frac{1+zf_0+z^2f_0+z^4tf_0-f_0}{z}=f_1={\frac {z \left( {z}^{2}+{z}^{3}f_0+1+zf_0-f_0 \right) 
	}{-{z}^{4}+{z}^{4}t-1+r_2z}}
\end{equation*}
from which $f_0$ can be computed:
\begin{align*}
f_0&= {\frac {-{z}^{4}-1+{z}^{4}t+2 z+\sqrt { z^8-2z^8t+2z^4+z^8t^2-2z^4t+1-4z^2}}{2z
		\left( -{z}^{4}t+1-2 z+{z}^{4} \right) }}\\
	&=1+z+2 z^2+3 z^3+(t+5) z^4+(8+2 t) z^5+(14+6 t) z^6+(23+12t)z^7+\cdots
\end{align*}
The special cases are as follows:
\begin{align*}
	f_0\big|_{t=0}&= {\frac {2 z+\sqrt { z^8+2z^4+1-4z^2}-{z}^{4}-1}{2z
			\left( +1-2 z+{z}^{4} \right) }}\\
	&=1+z+2z^2+3z^3+5z^4+8z^5+14z^6+23z^7+41z^8+69z^9+124z^{10}+\cdots
\end{align*}
which enumerates dispersed Dyck paths without $UUDD$ (sequence A191794) and
\begin{align*}
\frac{\partial f_0}{\partial t}\bigg|_{t=1}&=\frac{z^4}{(1-2z)\sqrt{1-4z^2}}\\
&=z^4+2 z^5+6 z^6+12 z^7+30 z^8+60 z^9+140 z^{10}+280 z^{11}+630 z^{12}+\cdots,
\end{align*}
which enumerates the number of occurrences of $UUDD$ in dispersed Dyck paths (sequence A100071).

Of course, dispersed Dyck paths with arbitrary endpoint can also be discussed using
$f_0+F(1)+G(1)+H(1)+K(1)$, but we are not going to display any formula; they are easy to obtain.

\bibliographystyle{plain}


\begin{thebibliography}{1}
	
	\bibitem{hackl.AofA}
	Benjamin Hackl, Clemens Heuberger, and Helmut Prodinger.
	\newblock {Counting Ascents in Generalized Dyck Paths}.
	\newblock In James~Allen Fill and Mark~Daniel Ward, editors, {\em 29th
		International Conference on Probabilistic, Combinatorial and Asymptotic
		Methods for the Analysis of Algorithms (AofA 2018)}, volume 110 of {\em
		Leibniz International Proceedings in Informatics (LIPIcs)}, pages
	26:1--26:15, Dagstuhl, Germany, 2018. Schloss Dagstuhl -- Leibniz-Zentrum
	f{\"u}r Informatik.
	
	\bibitem{theis}
	Kairi Kangro, Mozhgan Pourmoradnasseri, and Dirk~Oliver Theis.
	\newblock Short note on the number of 1-ascents in dispersed {D}yck paths.
	\newblock {\em Discrete Math. Algorithms Appl.}, 9(6):1750077, 8, 2017.
	
	\bibitem{Prodinger-kernel}
	Helmut Prodinger.
	\newblock The kernel method: a collection of examples.
	\newblock {\em S\'{e}m. Lothar. Combin.}, 50:Art. B50f, 19, 2003/04.
	
	\bibitem{Prodinger-prefix}
	Helmut Prodinger.
	\newblock Partial skew {D}yck paths: a kernel method approach.
	\newblock {\em Graphs Combin.}, 38(5):Paper No. 135, 11, 2022.
	
	\bibitem{Prodinger-JIS}
	Helmut Prodinger.
	\newblock Skew dyck paths having no peaks at level 1.
	\newblock {\em Journal of Integer Sequences}, 25(1):Article 22.1.6, 10 pages,
	2022.
	
	\bibitem{garden}
	Helmut Prodinger.
	\newblock A walk in my lattice path garden.
	\newblock {\em S\'{e}m. Lothar. Combin.}, 87b:49 p., 2023.
	
	\bibitem{OEIS}
	Neil J.~A. Sloane.
	\newblock The on-line encyclopedia of integer sequences.
	\newblock {\em Notices Amer. Math. Soc.}, 65(9):1062--1074, 2018.
	
\end{thebibliography}

\end{document}